\documentclass[12pt]{article}

\usepackage{amsmath}
\usepackage{amsfonts}
\usepackage{amsthm}
\usepackage{graphicx}
\usepackage{caption}
\usepackage{subcaption}
\usepackage{enumerate}
\usepackage{hyperref}
\usepackage[parfill]{parskip} 

\usepackage{fancyhdr}
\numberwithin{equation}{section}

\newcommand{\RR}{\mathbb{R}}      

\begin{document}

\title{Numerical and Symbolic Studies of the Peaceable Queens Problem}
\author{Yukun Yao \and Doron Zeilberger}

\maketitle

\begin{abstract}
We discuss, and make partial progress on, the peaceable queens problem, the protagonist of OEIS sequence {\tt A250000}. Symbolically, we prove that Jubin's construction of two pentagons is at least a local optimum. Numerically, we find the exact numerical optimums for some specific configurations. Our method can be easily applied to more complicated configurations with more parameters. 
\end{abstract}
\leavevmode
\\
\\
\\
{\bf Accompanying Maple package}

This  article is accompanied by a Maple package, {\tt PeaceableQueens.txt}, available from the url

{\tt http://sites.math.rutgers.edu/\~{}zeilberg/mamarim/mamarimhtml/peaceable.html}, 

where readers can also find an input and an output file, and nice pictures.

\section*{Introduction}

One of the fascinating problems described in the recent article [S2], about the great  {\bf On-Line Encyclopedia of Integer Sequences}, and in the beautiful and insightful video [S3]
is the {\it peaceable queens problem}.
It was chosen, by popular vote, to be assigned the {\bf milestone} `quarter-million' A-number, {\tt A250000}.

The question is the following:

{\it What is the maximal number, $m$, such that it is possible to place $m$ white queens and $m$ black queens on an
$n \times n$ chess board, so that no queen attacks a queen of the opposite color.}

Currently only thirteen terms are known:
\begin{align*}
n: && 1 && 2 && 3 && 4 && 5 && 6 && 7 && 8 && 9 && 10 && 11 && 12 && 13 \\
a(n): && 0 && 0 && 1 && 2 && 4 && 5 && 7 && 9 && 12 && 14 && 17 && 21 && 24
\end{align*}

In this paper, we'd like to consider this peaceable queens problem as a continuous question by normalizing the chess board to be the unit square $U:=[0,1]^2 = \{(x,y) \, | \, 0 \leq x, y \leq 1 \}$. Let $W \subseteq U$ be the region where white queens are located. Then the non-attacking region $B$ of $W$ can be defined as
$$
B = \{(x,y) \in U \,| \, \forall (u,v) \in W, x \neq u,  y \neq v, x+y \neq u+v,  y-x \neq v-u \}.
$$
So the continuous version of the peaceable queens problem is to find 
$$
\max_{W \in 2^U} ( \min(\textrm{Area} (W), \textrm{Area}(B)) ).
$$

Considering that the queen is able to move any number of squares vertically, horizontally and diagonally, it is reasonable to let $W$ be a convex polygon or a disjoint union of convex polygons whose boundary consists of vertical, horizontal and slope $\pm  1$ line segments, otherwise in many cases we can increase the area of white queens without decreasing the area of black queens. 

In this paper, we use a list $L$ of lists $[\,[a_1, b_1]\,,\, [a_2, b_2]\,,\, \dots\,,\, [a_n, b_n]\,]$ to denote the $n$-gon whose vertices are the $n$ pairs in the list $L$ and whose sides are the straight line segments connecting $[a_i, b_i]$ and $[a_{i+1}, b_{i+1}], (1 \leq i \leq n-1) $, and $[a_n, b_n]$ and $[a_1, b_1]$.

This paper is organized as follows. At first we look at Jubin's construction and prove that it is a local optimum. Though there is no rigorous proof, we conjecture and reasonably believe that it is indeed a global optimum at least for ``the continuous chess board", after numerous experiments with one, two and more components. Then we consider the optimal case under more restrictions, or under certain configurations, e.g., only one component or two identical squares or two identical triangles, etc. In some cases, the exact optimal parameters and areas can be obtained. Note that in this paper's figures, for convenience of demonstration, the color red is used to represent white queens and blue is for black queens.

\section*{Jubin's Construction}

As mentioned in [S2] (and [S1], sequence A250000), it is conjectured that Benoit Jubin's construction 
given in Fig. 5 of [S2], see also here:

{\tt http://sites.math.rutgers.edu/\~{}zeilberg/tokhniot/peaceable/P1.html}

or Figure 1, 
is  optimal for $n \geq 10$. Its value is $\lfloor \frac{7n^2}{48} \rfloor$.

\begin{figure}[h!]
  \includegraphics[width=\textwidth]{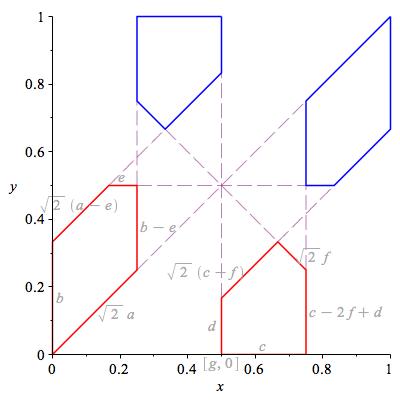}
  \caption{Benoit Jubin's Construction for a Unit Square}
\end{figure}

While we are, at present, unable to prove this,
we did manage to prove that when one generalizes Jubin's construction and replaces the
sides of the two pentagons with arbitrary parameters (of course subject to the obvious constraints so that both white and black queens reside in two pentagons),
then Jubin's construction is indeed (asymptotically) optimal, i.e. in the limit as $n$ goes to infinity.

{\bf Lemma 1} Normalizing the chess board to be  the unit square $\{(x,y) \, | \, 0 \leq x, y \leq 1 \}$,
if the white queens are placed in the union of the interiors of the following two pentagons
$$
[\, [0,0] \, , \, [a,a] \, , \, [a,a+b-e] \, , \, [a-e,a+b-e] \, , \, [0,b] \,   ],
$$
and
$$
[\, [g,0] \, , \, [g+c,0] \, , \, [g+c,c-2\,f+d] \, , \, [g+c-f,c-f+d] \, , \, [g,d] \,  ],
$$
where $a,b,c,d,e,f,g$ are between $0$ and $1$ and all coordinates and side lengths in Fig. 1 are non-negative and appropriate so that black queens also reside in two pentagons, then the black queens are located in the
interiors of the pentagons
$$
[\, [g,1] \, , \, [a,1] \, , \, [a,g+2\,c-2\,f+d-a] \, , \, 
[\frac{1}{2} \,g+ \frac{1}{2}\,d-\frac{1}{2}\,b+c-f, \frac{1}{2}\,g+ \frac{1}{2}\,d+ \frac{1}{2}\,b+c-f] \, , 
\,
$$
$$
[g,g+b] \, ] ,
$$
and
$$
[\, [1,1] \, , \, [g+c,g+c] \, , \, [g+c,a+b-e] \, , \, [a+b-e+g-d,a+b-e] \, , \, [1,1+d-g] \,  ]  .
$$

\begin{proof}
Since we only consider cases when the black queens also reside in two pentagons, this requirement provides natural constraints for these parameters $a,b,c,d,e,f,g$. Just to name a few, $a \leq g$ because the two pentagons are not overlapped, $g+c \leq 1$ because the right pentagon of white queens should entirely reside in the unit square, $d\leq g$ because otherwise the right pentagon of black queens will not exist and $c-f+d$, which is the $y$-coordinate of the highest point in the right pentagon of white queens, cannot be too large to ensure the right pentagon of black queens does not degenerate to a parallelogram. In these constraints we always use ``$\leq$" instead of ``$<$" so that the Lagrange multipliers will be able to work in a closed domain. 

With school geometry, it is obvious that black queens cannot reside in $0 \leq x < a$ since it is attacked by the left pentagon of white queens. Similar arguments work for the area $0<y \leq a+b$ and $g<x<g+c$. Now the leftover on the unit square is a union of two rectangles. By excluding $x+y < g+2c-2f+d, 0<y-x<b$ and $y-x<d-g$, these two rectangles are shaped into two pentagons and the coordinates of their vertices follow immediately.
\end{proof}

{\bf Lemma 2} The area of the white queens is
$$
ab \, - \,  \frac{1}{2}\,{e}^{2}+cd \, + \,  \frac{1}{2}\,{c}^{2}-{f}^{2} ,
$$
while the area of the black queens is
$$
-a-\frac{3}{4}\,{d}^{2}+2\,g-d-cd-ab-{f}^{2}- \frac{1}{2}\,{e}^{2}-\frac{3}{2}\,{c}^{2}+2\,bc-2\,af+3\,ac+2\,ad+2\,cf-ec-ed+be
$$
$$
+ae-bf+fd+ \frac{3}{2}\,bd-{a}^{2} - \frac{3}{4}\,{b}^{2}-2\,g c+ \frac{1}{2}\,g d- \frac{1}{2}\,g b+ag+gf- \frac{7}{4} \,{g}^{2} .
$$
\begin{proof}
For white queens, the left pentagon is a rectangle minus two triangles. Hence the area is 
$$
a(a+b-e) - \frac{1}{2} a^2 - \frac{1}{2} (a-e)^2 = ab - \frac{1}{2} e^2.
$$
The area of the right pentagon is 
$$
c(d+c-f) - \frac{1}{2} f^2 - \frac{1}{2} (c-f)^2 = \frac{1}{2} c^2 + cd - f^2.
$$
So the area of the white queens follows. 
For black queens, similarly, with the coordinates of the vertices in Lemma 1, simple calculation leads to the formula of its area. 
\end{proof}

{\bf Theorem} The optimal case of the two-pentagon configuration is Jubin's construction.
\begin{proof}

The procedure {\tt MaxC(L,v)} in the Maple package {\tt PeaceableQueens.txt} takes a list of length 2, $L$, consisting of polynomials in the list of variables $v$, and $v$ as inputs and outputs all the extreme points of $L[1]$, subject to the constraint $L[1]=L[2]$, using Lagrange multipliers.

Optimally, the areas of the white queens and black queens should be the same. Maximizing this quantity with the procedure {\tt MaxC} under this constraint shows that the maximum value is
$$
\frac{7}{48}  ,
$$
and this is indeed achieved by Jubin's construction, in which the white queens are located inside the pentagons
$$
[\, [0,0] \, , \, [\frac{1}{4}, \frac{1}{4}] \, , \,  [ \frac{1}{4}, \frac{1}{2}]
\, , \, [\frac{1}{6}, \frac{1}{2}] \, , \, [0,\frac{1}{3}] \,  ] ,
$$
and
$$
[\, [ \frac{1}{2},0] \, , \, [\frac{3}{4},0] \, , \, [\frac{3}{4},\frac{1}{4}] \, , \,
[\frac{2}{3},\frac{1}{3}] \, , \, [ \frac{1}{2},\frac{1}{6}] \, ] ,
$$
and the black queens reside inside the pentagons
$$
[\, [ \frac{1}{2},1] \, , \, [\frac{1}{4},1]
\, , \, [\frac{1}{4},\frac{3}{4}] \, , \, [\frac{1}{3},\frac{2}{3}] \, , \,
[ \frac{1}{2},\frac{5}{6}] \,  ] ,
$$
and
$$
[\, [1,1] \, ,\ \, [\frac{3}{4},\frac{3}{4}] \, , \, [\frac{3}{4}, \frac{1}{2}] \, , \,
[\frac{5}{6}, \frac{1}{2}] \, , \, [1,\frac{2}{3}] \,  ] . 
$$
\end{proof}

It seems natural that two components are optimal because if there is only one connected component for white queens, black queens still have two connected components. From the view of symmetry, it seems good to add the other component for white queens. In the rest of the paper, it is shown that with only one connected component it is unlikely to surpass the $\frac{7}{48}$ result. And by experimenting with three or more connected components for the white queens, it seems that it is not possible to improve on
Jubin's construction, hence we believe that it is indeed optimal (at least asymptotically). By the way, Donald Knuth kindly informed us that
what we (and the OEIS) call Jubin's construction already appears in  Stephen Ainley's delightful book "Mathematical Puzzles"[A], p. 31, Fig, 28(A) .

\section*{Single Connected Component}
In this section, we try to find the optimal case when $W$ is a single connected component and when the configuration is restricted to rectangles, parallelograms, triangles and finally obtain a lower bound for the optimal case of one connected component.

\subsection*{A Single Rectangle}
Let the rectangle for white queens be $[\,[0,0]\,, \,[a,0]\,, \,[a,b]\,,\, [0,b]\,]$, with the obvious fact that for a rectangle with a given size, placing it in the corner will lead to the largest non-attacking area. The area for white queens is $ab$ and the area for black queens is $(1-a-b)^2$. We'd like to find the maximum of $ab$ under the condition 
$$
ab = (\max(1-a-b, 0))^2, \quad 0 \leq a,b \leq 1.
$$

\begin{figure}[h!]
  \includegraphics[width=\textwidth]{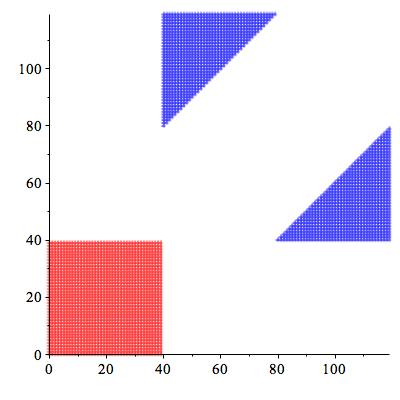}
  \caption{The Optimal Rectangle for a 120 by 120 Chess Board}
\end{figure}

Since $a$ and $b$ are symmetric, the maximum must be on the line $a=b$. Hence the optimal case is when 
$$
a=b=\frac{1}{3}
$$
and the largest area for peaceable queens when the configuration for white queens is a rectangle is $\frac{1}{9}$.

\subsection*{A Single Parallelogram}
Let the parallelogram for white queens be $[\,[0,0]\,,\,[a,a]\,,\,[a,a+b]\,,\,[0,b]\,]$.

\begin{figure}[h!]
  \includegraphics[width=\textwidth]{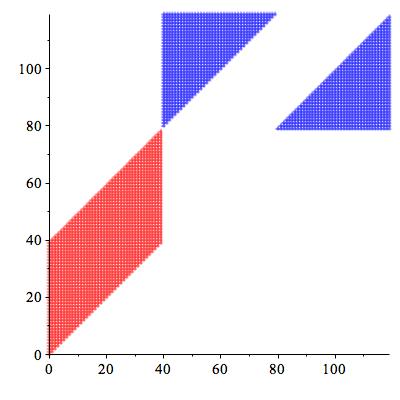}
  \caption{The Optimal Parallelogram for a 120 by 120 Chess Board}
\end{figure}

Note that as mentioned in the beginning of this section, because the line segment must be vertical, horizontal or of slope $\pm 1$ and the corner is the best place to locate a shape, there are only two kinds of parallelograms, the other one being $[[0,0],[b,0],[a+b,a],[a,a]]$. Obviously they are symmetric with respect to the line $y=x$, so let's focus on one of them. 

The area for white queens is still $ab$ and the area for black queens is still $(\max(1-a-b, 0))^2$. So similarly with the rectangle case, the maximum area $\frac{1}{9}$ is reached when
$$
a=b=\frac{1}{3}.
$$

\subsection*{A Single Triangle}
With similar arguments as in the last subsection, the optimal triangle must have the format: $[\,[0,0]\,,\,[0,a]\,,\,[a,a]\,]$. The area for white queens is 
$$
\frac{1}{2} a^2
$$
and the area for black queens is 
$$
\frac{1}{2} (1-a)^2.
$$
with the condition $0 \leq a \leq 1$.

Hence, when $a=\frac{1}{2}$ the area reaches its maximum $\frac{1}{8}$, which is better than the rectangle or parallelogram configuration. 

\begin{figure}[h!]
  \includegraphics[width=\textwidth]{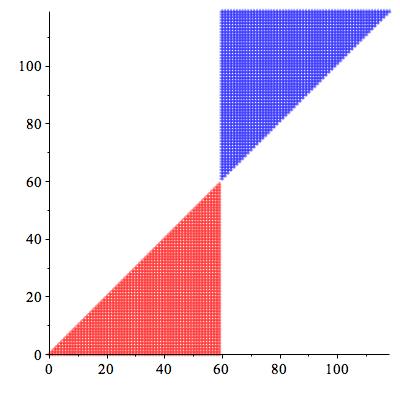}
  \caption{The Optimal Triangle for a 120 by 120 Chess Board}
\end{figure}

By the way, $[\,[0,0]\,,\,[0,a]\,,\,[a,0]\,]$ won't be a good candidate for optimal triangles because we can always extend it to a square $[\,[0,0]\,,\,[0,a]\,,\, [a,a] \, , \, [a,0]\,]$ without decreasing the area of black queens. Then its maximum cannot exceed the maximum of rectangles, which is $\frac{1}{9}$.

\subsection*{A Single Hexagon}
After looking at specific configurations in the above three subsections, we'd like to find some numerical lower bounds for the single connected component configuration. It is interesting to find out or at least get a numerical estimation how large the area of white or black queens can be if the white queens are in a single connected component.  Note that from rectangles and parallelogram we get a lower bound $\frac{1}{9} \approx 0.1111$ and from triangles we get a better lower bound $\frac{1}{8} = 0.125$.

The natural thing is that we want to place the polygon in a corner. Because of the restriction of the orientations of its sides, at most it can be an octagon. Let's place the polygon in the lower left corner. Then we immediately realize that it is a waste if the polygon doesn't fill the lower left corner of the unit square. It is the same for the upper right side of the polygon. If part of its vertices are $[\,[a,b]\,,\, [a, b+c]\,,\,[a-d, b+c+d]\,,\,[a-d-f, b+c+d]\,]$, then we can always extend the polygon to $[\, \dots \,,\,[a,b]\,,\, [a, b+c+d]\,,\, [a-d-f, b+c+d]\,, \, \dots \, ] $ without decreasing the area of black queens. 

Hence the general shape is a hexagon 
$$
[\, [0,0]\,,\,[a,0]\,,\,[a+b,b]\,,\,[a+b, b+c]\,,\,[d, b+c]\,,\,[0,b+c-d]\,]
$$ 
with four parameters. Then the area for white queens is 
$$
(a+b)(b+c) - \frac{1}{2} (b^2+d^2),
$$
and the area for black queens is 
$$
\frac{1}{2} (1-a-b-c)^2 + \frac{1}{2} (1-a-2b-c+d)^2.
$$
With the procedure {\tt MaxC}, one of the local maximums found using Lagrange multipliers is when 
$$
a=c=d = \frac{1}{2}, \, b=0.
$$
However, actually this is the optimal triangle with an area of $\frac{1}{8}$.

\begin{figure}[h!]
  \includegraphics[width=\textwidth]{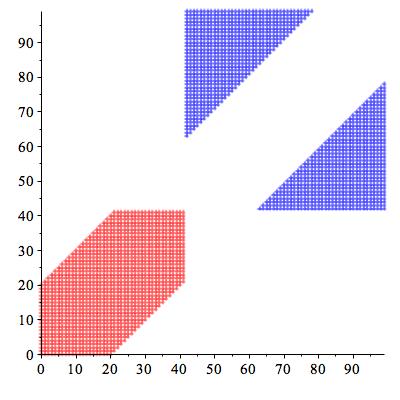}
  \caption{The Nearly Best Lower Bound Configuration for a 100 by 100 Chess Board}
\end{figure}

Another local maximum is when $a=b=c=d$. In that case, we have 
$$
3a^2 = (1-3a)^2.
$$
Hence when 
$$a = \frac{3-\sqrt{3}}{6} \approx 0.2113248654,$$ the area of white queens is maximized at 
$$ 3a^2 = \frac{2-\sqrt{3}}{2} \approx 0.1339745962. $$
The best configuration of hexagons is found and at least we have a numerical lower bound 0.1339745962 for the best single component configuration.

\section*{Two Components}
Since in Jubin's construction, there are two pentagons, it is natural to think of the optimum of certain two-component configurations. The difficulty for analyzing the two-component is that more parameters are introduced and the area formula for black queens becomes a much more complicated piece-wise function. 

In this section, the cylindrical algebraic decomposition algorithm in quantifier elimination is applied to find out the exact optimal parameters and the maximum areas. Given a set $S$ of polynomials in $\RR^n$, a cylindrical algebraic decomposition is a decomposition of $\RR^n$ into semi-algebraic connected sets called cells, on which each polynomial has constant sign, either +, - or 0. With such a decomposition it is easy to give a solution of a system of inequalities and equations defined by the polynomials, i.e. a real polynomial system.

\subsection*{Two Identical Squares}
To keep the number of parameters as few as possible, the configuration of two identical squares is the first we'd like to study. There are two parameters, the side length $a$ and the $x$-coordinate $s$ of the lower left vertex of the right square, the left square's lower left vertex being the origin. 

The two squares are 
$$[\, [0,0] \,,\, [a,0] \,,\, [a,a] \,,\, [0,a] \, ]$$
and 
$$[\, [s,0] \,,\, [s+a,0] \,,\, [s+a,a] \,,\, [s,a] \, ].$$ 

\begin{figure}[h!]
  \includegraphics[width=\textwidth]{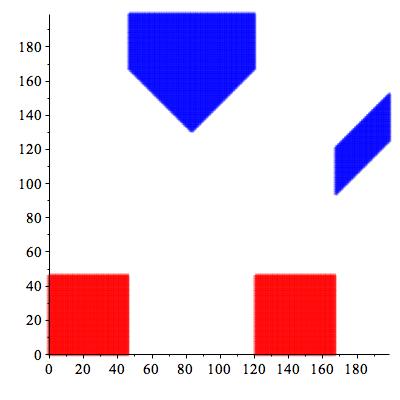}
  \caption{The Nearly Optimal Two Identical Squares Configuration for a 200 by 200 Chess Board}
\end{figure}

Based on this configuration, the domain is 
$$
0 \leq a \leq \frac{1}{2}, \quad a \leq s \leq 1-a.
$$

The area of white queens is 
$$
2a^2.
$$

Actually the formula for black queens is very complicated, especially when $a$ is small there may be a lot of components for $B$. However, by experimentation (procedure {\tt FindM2Square}), we found that for all mid-range $s \in [0.24, 0.76]$, $a$ around 0.23 will always maximize the area. Then we just need to focus on the shape of $B$ when $a$ is not far from its optimum. 

The area of black queens is 
$$
(s-a)(1-s-a) + \frac{1}{4}(s-a)^2 + (\max(1-s-2a, 0))^2 + \max(s-2a,0)(1-s-a).
$$

The domain for $a$ and $s$ is a triangle. The area formula for black queens shows that the two lines $s=2a$ and $s=1-2a$ separate the domain into 4 regions. In each region, we have a polynomial formula for the area of black queens. Since the area of white queens $W$ is just a simple formula of $a$, we need to maximize $a$ with the condition $W=B$.

When $s \geq 2a$ and $s \geq 1-2a$, by cylindrical algebraic decomposition we obtained 
$$
\begin{cases} 
    \frac{1}{2}(-1+\sqrt{2}) \leq a < \frac{1}{27}(1+2\sqrt{7}) & s=\frac{4+a}{7} + \frac{2}{7} \sqrt{4-19a+9a^2} \\
    \frac{1}{27}(1+2\sqrt{7}) \leq a < \frac{1}{18}(19-\sqrt{217}) &s=\frac{4+a}{7} \pm \frac{2}{7} \sqrt{4-19a+9a^2} \\
    a=\frac{1}{18}(19-\sqrt{217}) & s= \frac{4+a}{7} - \frac{2}{7} \sqrt{4-19a+9a^2}
\end{cases}.
$$

When $s \leq 2a$ and $s \geq 1-2a$, the result is an empty set. 

When $s \leq 2a$ and $s \leq 1-2a$, we obtained
$$
\frac{2}{9} \leq a \leq \frac{1}{7}(3-\sqrt{2}), \quad s=2-7a-2\sqrt{-2a+9a^2}.
$$

When $s \geq 2a$ and $s \leq 1-2a$, we obtained
$$
\frac{2}{9} \leq a \leq \frac{1}{27}(1+2\sqrt{7}), \quad s=3a - \frac{2}{\sqrt{3}} \sqrt{1-7a+12a^2}.
$$

Comparing the four cases, we found that the largest area occurred in case 1, when 
$$
a = \frac{1}{18} (19-\sqrt{217}) \approx 0.2371711193,
$$
$$
s = \frac{13}{18} - \frac{1}{126} \sqrt{217} \approx 0.6053101598 .
$$
The largest area is ${\frac {289}{81}}-{\frac {19\,\sqrt {217}}{81}} \approx  0.112500281.$

\subsection*{Two Identical Triangles}
The configuration of two identical isosceles right triangles with the same orientation is the next to be considered. There are also two parameters, the leg length $a$ and the $x$-coordinate $s$ of the lower left vertex of the triangle on the right. Note that the slopes of both triangles' hypotenuses are $+1$. 

The two isosceles right triangles are
$$[\,[0,0]\,,\,[a,0]\,,\,[a,a]\,]$$
and 
$$
[\,[s,0]\,,\,[a+s,0]\,,\,[a+s,a]\,].
$$
The domain for the two parameters $a$ and $s$ is also 
$$
0 \leq a \leq \frac{1}{2}, \quad a \leq s \leq 1-a.
$$
The area of white queens is $a^2$ and for the area of black queens, by numerical experimentation, we found that for all mid-range $s \in [0.32, 0.68]$, the area is maximized when $a$ is around 0.31. Hence for $a$ around 0.31, we have that the area of black queens is 
$$
2(s-a)(1-s-a) + \frac{1}{2} (s-a)^2 + \frac{1}{2} (1-s-a)^2 + \frac{1}{2} (\max (1-s-2a, 0))^2.
$$
When $s \geq 1-2a$, by cylindrical algebraic decomposition we obtained 
$$
\begin{cases} 
    \frac{1}{2}(2-\sqrt{2}) \leq a < \frac{1}{4}(-1+\sqrt{5}) & s=\frac{1}{2} + \frac{1}{2} \sqrt{3-12a+8a^2} \\
   \frac{1}{4}(-1+\sqrt{5}) \leq a < \frac{1}{4}(3-\sqrt{3}) &s=\frac{1}{2} \pm \frac{1}{2} \sqrt{3-12a+8a^2} \\
    a=\frac{1}{4}(3-\sqrt{3}) & s=\frac{1}{2} - \frac{1}{2} \sqrt{3-12a+8a^2}
\end{cases}.
$$
When $s \leq 1-2a$, we obtained 
$$
\frac{1}{11}(5-\sqrt{3}) \leq a \leq \frac{1}{4}(-1+\sqrt{5}), s=2a - \sqrt{2-10a+12a^2}.
$$
Hence the area is maximized when
$$a = \frac{1}{4}(3-\sqrt{3}) \approx 0.316987298,$$
$$
s = \frac{1}{2}.
$$
The largest area is $\frac{3}{4} - \frac{3}{8} \sqrt{3} \approx 0.1004809470$.

\begin{figure}[h!]
  \includegraphics[width=\textwidth]{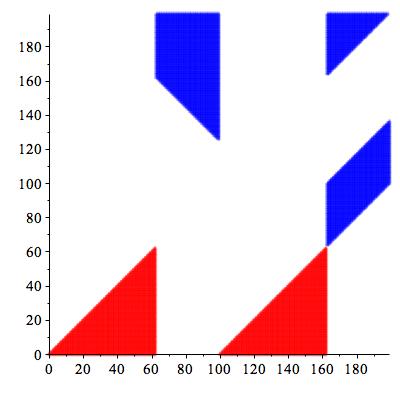}
  \caption{The Nearly Optimal Two Identical Isosceles Right Triangles with the Same Orientation Configuration for a 200 by 200 Chess Board}
\end{figure}

Thanks to the referee's suggestions, a larger area can be obtained if two identical isosceles right triangles with different orientations are considered. For example, if we take the two triangles to be $$[\,[0,0]\,,\,[a,0]\,,\,[a,a]\,]$$
and 
$$
[\,[1-a,0]\,,\,[1,0]\,,\,[1-a,a]\,],
$$
then the area of black queens is 
$$
a(1-2a) + (\frac{1}{2} - a)^2 = -a^2 + \frac{1}{4}.
$$
Equalizing the areas of white queens and black queens, we get 
$$
\textrm{Area}(W) = a^2 = \frac{1}{8},
$$
which is greater than the optimal case of two identical isosceles right triangles with the same orientation.

\begin{figure}[h!]
  \includegraphics[width=\textwidth]{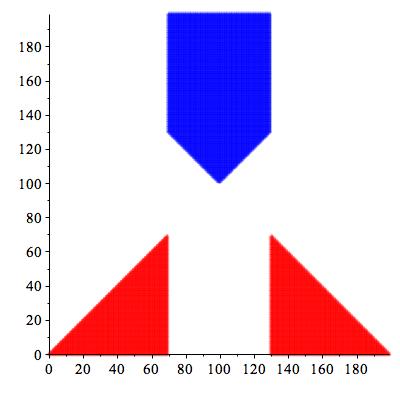}
  \caption{An Example of Two Identical Isosceles Right Triangles with Different Orientations Configuration for a 200 by 200 Chess Board}
\end{figure}

\subsection*{One Square and One Triangle with the Same Side Length}

\begin{figure}[h!]
  \includegraphics[width=\textwidth]{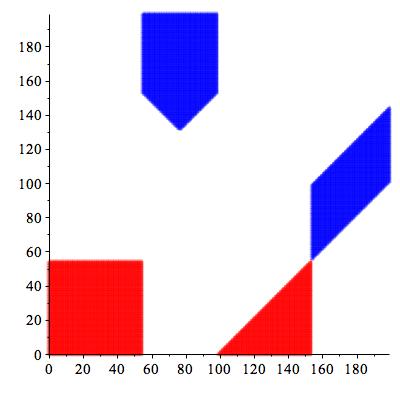}
  \caption{The Nearly Optimal One Square and One Triangle (with the same side length) Configuration for a 200 by 200 Chess Board}
\end{figure}

With the same notations as the above two subsections, let $W$ be the union of the square 
$$[\, [0,0] \,,\, [a,0] \,,\, [a,a] \,,\, [0,a] \, ]$$
and the triangle
$$ [\,[s,0]\,,\,[a+s,0]\,,\,[a+s,a]\,]. $$
Then the area of white queens is $\frac{3}{2} a^2$ and the area of black queens is 
$$
a(s-a)(1-s-a) + \frac{1}{4} (s-a)^2 + (\max(1-s-2a, 0))^2
$$
when $a$ is around its optimum 0.27 and $s \in [0.28, 0.72]$. It is obtained that when $s \geq 1-2a$
$$
\begin{cases} 
    \frac{1}{2}(-2+\sqrt{6}) \leq a < \frac{1}{21}(1+\sqrt{22}) & s=\frac{4-a}{7} + \frac{1}{7} \sqrt{16-64a+22a^2} \\
   \frac{1}{21}(1+\sqrt{22}) \leq a < \frac{2}{11}(8-\sqrt{42}) &s=\frac{4-a}{7} \pm \frac{1}{7} \sqrt{16-64a+22a^2} \\
    a=\frac{2}{11}(8-\sqrt{42}) & s=\frac{4-a}{7} - \frac{1}{7} \sqrt{16-64a+22a^2}
\end{cases},
$$
and when $s \leq 1-2a$
$$
\frac{1}{15} (6-\sqrt{6}) \leq a \leq \frac{1}{21} (1+\sqrt{22}), \quad s = \frac{7a}{3} - \frac{1}{3} \sqrt{12-72a+106a^2}.
$$
Consequently, we have the maximized area when
$$
a = \frac{2}{11} (8-\sqrt{42}) \approx 0.276228965，
$$
$$
s = \frac{112}{33} - \frac{14}{33} \sqrt{42} - \frac{50}{33} \sqrt{7} + \frac{52}{33} \sqrt{6} \approx 0.495622162.
$$
The largest area is $\frac{636}{121} - \frac{96}{121} \sqrt{42} \approx 0.1144536616$. Among the three configurations in this section, we found that this configuration with one square and one triangle has the largest area. 

\section*{Future Work and Final Remarks}
Our method can be easily generalized for configurations with more components and/or more parameters. For instance, let's consider the configuration of two squares, not necessarily identical. Then there are three parameters, the side length $a$ of the left square, the side length $b$ of the right square and the $x$-coordinate $s$ of the right square's lower left vertex. For fixed $b$ and $s$, we can find the interval of $a$ in which the optimum is located. Then for each fixed $s$, we are able to find the interval of $b$ such that its corresponding $a$ will lead to the largest area $a^2 + b^2$. When the estimated optimal parameters are determined, a piece-wise function of the area of black queens follows. 

The main difficulty of this peaceable queens problem lies in the number of parameters and the complexity of the area formula of black queens. When there are multiple components, as long as the number of parameters is limited, it should be still doable. For example, the configuration of three identical squares which are placed equidistantly has only one parameter, the side length $a$. When the chess board is 240 by 240, the optimal $a$ is around 40, which means in the unit square the optimal side length is around $\frac{1}{6}.$ 

In conclusion, in this paper we prove that Jubin's configuration is a local optimum. Optimal cases of some certain configurations are discussed. Future work includes the exact solution of complicated configurations with numerous parameters, whether the white queens have two components under the best configuration, and proof or disproof that Jubin's configuration is indeed the best.

\bigskip
\bigskip
\bigskip
\bigskip

{\bf Acknowledgement} We are thankful to Neil J.A. Sloane for introducing us to this interesting problem and we are also grateful to Lun Zhang for pleasant conversations and helpful remarks. We sincerely thank the anonymous referee for careful review and useful suggestions. We appreciate Shalosh B. Ekhad's impeccable computing support. 

\bigskip
\bigskip
\bigskip
\bigskip
\bigskip
\bigskip
\bigskip

{\bf References}

[A] Stephen Ainley, {\it Matematical puzzles}, Prentice Hall, Upper Saddle River, NJ, 1977.

[BPR] Saugata Basu, Richard Pollack and  Marie-Françoise Roy, {\it Algorithms in real algebraic geometry}, second edition. Algorithms and Computation in Mathematics {\bf 10}, Springer-Verlag, Berlin, 2006.

[Bo] Robert Bosch, {\it Peaceably coexisting armies of queens}, Optima (Newsletter of the Mathematical Programming Society) {\bf 62.6-9}: 271, 1999.

[Br] Christopher W. Brown, {\it Simple CAD Construction and Its Applications},  Journal of Symbolic Computation, {\bf 31}, 521-547, 2001.

[CJ] Bob Caviness and Jeremy Johnson (Eds.), {\it Quantifier Elimination and Cylindrical Algebraic Decomposition}, Springer-Verlag, New York, 1998.

[J] Benoit Jubin, {\it Improved lower bound for A250000}, \hfill\break
{\tt https://oeis.org/A250000/a250000\_1.txt}

[K] Donald E. Knuth, {\it Satisfiability, Fascicle 6, volume 4 of The Art of Computer Programming}, Addison-Wesley, 2015.

[S1] Neil J.A. Sloane, {\it The On-Line Encyclopedia of Integer Sequences}, \hfill\break
{\tt https://oeis.org/}  .

[S2] Neil J.A. Sloane, {\it The On-Line Encyclopedia of Integer Sequences}, Notices of the American Mathematical
Society {\bf 65} \#9 , 1062-1074, 2018.

[S3] Neil J.A. Sloane, {\it Peaceable Queens - Numberphile} video, available from \hfill\break
{\tt https://www.youtube.com/watch?v=IN1fPtY9jYg}

\bigskip
\hrule
\bigskip

Yukun Yao, Department of Mathematics, Rutgers University (New Brunswick), Hill Center-Busch Campus, 110 Frelinghuysen
Rd., Piscataway, NJ 08854-8019, USA. \hfill\break
Email: {\tt yao at math dot rutgers dot edu}    .
\bigskip

Doron Zeilberger, Department of Mathematics, Rutgers University (New Brunswick), Hill Center-Busch Campus, 110 Frelinghuysen
Rd., Piscataway, NJ 08854-8019, USA. \hfill\break
Email: {\tt DoronZeil at gmail dot com}    .

\end{document}